\documentclass[11pt]{article}

\usepackage{arxiv}
\usepackage[utf8]{inputenc}
\usepackage[T1]{fontenc}
\usepackage{lmodern}
\usepackage{amsmath,amssymb,amsfonts,mathtools}
\usepackage{booktabs}
\usepackage{array}
\usepackage{graphicx}
\usepackage{natbib}
\usepackage{microtype}
\usepackage{enumitem}
\usepackage{listings}
\usepackage{caption}
\usepackage{subcaption}

\setlist[itemize]{leftmargin=1.5em}
\setlist[enumerate]{leftmargin=1.8em}
\newcommand{\creg}{C_{\mathrm{reg}}}

\newcommand{\R}{\mathbb{R}}
\newcommand{\D}{\mathcal{D}}

\title{Diagnostic Certificates of Data Quality and Regression Identifiability for Koopman Identification}
\author[1,2]{Yue Wu}
\affil[1]{School of Automation Science and Engineering, Xi'an Jiaotong University, Xi'an 710049, China}
\affil[2]{Xinjiang Cigarette Factory, Hongyun Honghe Tobacco (Group) Co., Ltd., Urumqi 830000, China}
\date{\today}

\begin{document}
\maketitle

\begin{abstract}
Classical persistent excitation criteria usually assess whether an input or regressor signal is sufficiently rich. In Koopman and EDMD with control (EDMDc), however, data quality is determined by the concatenation of lifted state features and control inputs. Input-rich data can still visit a narrow state region, well-spread state samples can still produce degenerate lifted features, and both can fail to condition the final regression problem.

This paper develops a diagnostic certificate framework for locating these failures. The certificates separate state-space coverage and clustering, lifted-feature nondegeneracy, and the final regression spectrum. The regression-spectrum certificate is the layer with direct theoretical guarantees: it controls the active standardized design's smallest singular value, has Fisher-information and one-step EDMDc stability interpretations, and admits a finite-sample lower bound under a population spectral gap. We also give structural examples and a Schur-complement condition showing why state, lifted, input, and regression diagnostics cannot be substituted for one another.

As a sampling example, IGPE-DOPT uses these certificates to score candidate trajectory segments. Experiments on Duffing, Van der Pol, and Lorenz systems compare input-, state-, lifted-, and regression-oriented baselines. The results show that certificate layers separate, budget and weights shift bottlenecks, and downstream prediction or control performance is not monotone in any single certificate. The framework is therefore intended as an interpretable diagnostic and data-collection guide, not as a universal optimality claim.
\end{abstract}

\keywords{Koopman operator; EDMDc; data quality; persistent excitation; geometric certificates; regression identifiability; optimal experimental design}

\section{Introduction}

Nonlinear dynamical systems arise broadly in fluids, chemical processes, robotics, energy systems, and biological systems. Prediction and control of such systems require a balance between model expressiveness and computational tractability. The Koopman operator framework introduces tools from linear systems theory into nonlinear dynamics by lifting nonlinear state evolution to a linear operator acting on observable functions. Extended dynamic mode decomposition (EDMD) and EDMD with control (EDMDc) turn this idea into finite-dimensional regression problems.

For controlled systems, EDMDc fits a linear map from the current lifted state and the control input to the next lifted state. The quality of the identified model is therefore governed by a joint regression object: the visited states, the dictionary evaluations on those states, and the applied controls all enter the same least-squares design. This point is easy to obscure if classical persistent-excitation language is imported too directly. A signal can be rich in the input channel while the realized state trajectory remains confined to a narrow region. State samples can appear well spread while the selected dictionary is nearly constant or highly correlated on that region. Lifted features can be individually nondegenerate while the final lifted-state/control regression remains ill conditioned because the control columns are redundant with the lifted-state columns.

This structure reveals a multilayer data-quality problem in Koopman identification. The first layer is state-space geometry: whether trajectories cover the target region, whether directions are sufficiently diverse, and whether samples cluster locally. The second layer is lifted-feature geometry: whether the fixed dictionary forms a nondegenerate active feature space on these state samples. The third layer is regression-space geometry: whether lifted features and control inputs, after concatenation, form a regression matrix that can be estimated stably. These layers are related, but they are not interchangeable.

The goal of this paper is not to redefine classical persistent excitation or to propose a universally optimal sampler. The goal is more specific: for data collection in Koopman/EDMDc with a fixed dictionary, we provide computable and ablatable diagnostic certificates that separate where data quality is gained or lost. The paper emphasizes two principles. First, the spectral condition of the final active standardized regression matrix should be inspected directly, because that is the layer closest to estimation stability. Second, upstream geometric certificates should still be reported, because they explain why the final regression object becomes ill conditioned and which layer should be repaired by subsequent data collection.

The empirical part is therefore organized as a diagnostic evidence chain rather than as a winner-takes-all benchmark. We compare IGPE-DOPT with classical input-probing and design baselines, including adaptive PE-style probing (A-PE), deterministic multisine/chirp optimal-input probing (OID), state-geometric GPE collection, random and Sobol sampling, state-space k-center selection, lifted D-optimal design, regression D-optimal design, and regression E-optimal design. These baselines test whether the proposed certificates explain the differences among input-rich, state-covering, lifted-space, and regression-space acquisition rules.

The contributions are as follows.
\begin{enumerate}
  \item \textbf{Multilayer diagnostic certificates.} We define directional coverage, Frostman-type non-clustering, state isotropy, lifted-feature isotropy, regression isotropy, and a bottleneck-style composite certificate. We specify their standardization, active-column handling, and interpretation boundaries.
  \item \textbf{Regression-identifiability theory.} We prove that the final regression-spectrum certificate is directly connected to the minimum singular value of the active standardized design, a Fisher-type information quantity, and one-step EDMDc estimation stability under a fixed design matrix. We also give a population-to-sample lower bound for the regression spectrum.
  \item \textbf{Non-interchangeability of layers.} We use counterexamples to show that state geometry does not imply lifted or regression geometry, and we use a Schur-complement condition to show that nondegenerate lifted features or input variance alone is insufficient to guarantee nondegeneracy of the full regression matrix.
  \item \textbf{Certificate-driven sampling and an evidence-chain experiment.} We use IGPE-DOPT to show how certificates enter data acquisition. Using the real project evidence in Tables~\ref{tab:cert-hierarchy}--\ref{tab:weight-sensitivity} and Figures~\ref{fig:budget}--\ref{fig:a1}, we demonstrate certificate-layer separation, regression-theory correspondence, budget and weight sensitivity, and nonmonotonicity in downstream tasks.
\end{enumerate}

The rest of the paper is organized as follows. Section~\ref{sec:related} reviews persistent excitation, data informativity, space-filling design, and Koopman data acquisition. Section~\ref{sec:edmdc} gives the standardized EDMDc formulation. Section~\ref{sec:certificates} defines the diagnostic certificates. Section~\ref{sec:theory} gives theoretical guarantees and boundaries. Section~\ref{sec:igpe} introduces IGPE-DOPT and the comparison methods. Section~\ref{sec:experiments} reports experiments. Section~\ref{sec:discussion} discusses limitations. Section~\ref{sec:conclusion} concludes.

\section{Related Work}
\label{sec:related}

\subsection{Persistent excitation and data informativity}

Classical PE theory shows that, in adaptive control and linear system identification, parameter directions can be sufficiently excited if the input or regressor signal satisfies excitation conditions of sufficient order \citep{narendra1987persistent,willems2005note}. The data-informativity framework further emphasizes that whether data are sufficient depends on the target task: stabilization, prediction, and full model identification require different information \citep{van2020data}. This paper adopts this task-related viewpoint, but specializes the task object to the active standardized regression matrix of EDMDc.

\subsection{Information matrices and optimal experimental design}

Optimal experimental design has long used Fisher information matrices, D-optimality, and E-optimality to characterize parameter-estimation quality \citep{hjalmarsson2005experiment,deflorian2011design,decock2016doptimal,wilson2014trajectory}. These methods are, at their core, spectral-geometric methods for regression matrices. The certificate $\creg$ is consistent with this tradition, but it is embedded in the lifted-feature/input concatenation structure of EDMDc. The distinction in this paper is the explicit separation of state, lifted-feature, and regression diagnostic layers, rather than reducing all diagnostics to a single information-matrix score.

The experimental baselines reflect this lineage. The OID baseline uses deterministic multisine/chirp probing as a practical optimal-input-design surrogate. The REG-DOPT and REG-EOPT baselines isolate log-determinant and minimum-eigenvalue objectives in the final regression space, while LIFT-DOPT applies the same design logic before input concatenation. These are not treated as straw-man baselines: they represent plausible ways a practitioner might collect data if the only objective were input richness, lifted-feature diversity, or regression information.

\subsection{Space-filling design and state coverage}

Space-filling input design for nonlinear system identification shows that a rich input spectrum is insufficient to guarantee nonlinear model quality; coverage of the state space or of the joint input-state space is also important \citep{kiss2024spacefilling,liu2025spacefilling,smits2024spacefilling}. We incorporate this viewpoint but emphasize that state coverage is only an upstream diagnostic. For EDMDc, dictionary lifting and concatenated regression can reintroduce degeneracy, so lifted-feature geometry and regression geometry must still be inspected.

\subsection{Koopman/EDMDc identification and data acquisition}

The Koopman/EDMD literature has studied active learning, generalized PE, fundamental-lemma-type results with Koopman embeddings, regularization, and sampling distributions \citep{williams2015datadriven,proctor2016dynamic,korda2018linear,brunton2022modern,abraham2019active,boddupalli2019koopman,shang2024willems,dahdah2022system,philipp2025error}. These works provide different views of data requirements in Koopman identification. This paper focuses on reportable diagnostics for a finite dataset: given the data, dictionary, and inputs, how can one determine whether the data-quality bottleneck lies in the state space, lifted space, or final regression space?

Our GPE-STATE baseline is closest in spirit to geometry-aware active exploration: it uses state-space certificate improvement as the main acquisition signal. In contrast, IGPE-DOPT mixes state, lifted, and regression objectives. This distinction is important for positioning the algorithmic contribution. IGPE-DOPT is not presented as a universal active-learning replacement for prior Koopman acquisition methods. It is an instantiation showing how the proposed diagnostics can be turned into a data-collection rule and compared against input-, state-, lifted-, and regression-oriented alternatives.

\subsection{Coverage exploration, finite excitation, and control tasks}

Results on ergodic exploration, coverage control, finite excitation, and fast stabilization in reinforcement learning all show that data quality can be characterized through coverage, information-matrix spectra, or minimum eigenvalues of regression matrices \citep{miller2016ergodic,rickenbach2024active,parikh2019integral,lale2022reinforcement}. We apply these ideas to hierarchical diagnostics for Koopman/EDMDc, but we do not compress all certificates into a new PE condition.

\section{Problem Formulation and Standardized EDMDc Framework}
\label{sec:edmdc}

Consider a discrete-time controlled system
\begin{equation}
x_{k+1}=f(x_k,u_k), \qquad
x_k\in\R^{n_x}, \quad u_k\in\R^{n_u}.
\end{equation}
The dataset is
\begin{equation}
\D_N=\{(x_k,u_k,x_{k+1})\}_{k=1}^{N}.
\end{equation}
For a dictionary $\psi:\R^{n_x}\to\R^{d_\psi}$, define the augmented regressor
\begin{equation}
\xi_k=
\begin{bmatrix}
\psi(x_k)\\
u_k
\end{bmatrix}
\in\R^{p},
\qquad
p=d_\psi+n_u.
\end{equation}
Stacking samples gives
\begin{equation}
\Phi=
\begin{bmatrix}
\xi_1^\top\\
\cdots\\
\xi_N^\top
\end{bmatrix},
\qquad
Y=
\begin{bmatrix}
\psi(x_2)^\top\\
\cdots\\
\psi(x_{N+1})^\top
\end{bmatrix}.
\end{equation}
EDMDc solves
\begin{equation}
\min_K \|Y-\Phi K\|_F^2.
\end{equation}
The geometry of this problem is determined by $\Phi^\top\Phi$. If $\Phi^\top\Phi$ has small eigenvalues, some lifted/input directions have not been sufficiently excited by the data, and the least-squares estimate is sensitive to noise. Classical input PE inspects only the input sequence and does not directly guarantee the spectral condition of this concatenated matrix.

\subsection{Active standardized coordinates}

The columns of the raw regression matrix $\Phi$ may have very different scales. For example, polynomial features can be much larger than input columns, and a constant observable produces a zero-variance column. This paper defines spectral certificates in centered and active-column standardized coordinates. Write
\begin{equation}
\Phi=\mathbf{1}\mu_\Phi^\top+\bar{\Phi}D_\Phi,
\end{equation}
where $\mu_\Phi$ is the column mean, $D_\Phi$ is the active-column scale matrix, and $\bar{\Phi}$ is the active standardized regression matrix. If a column variance is below a numerical threshold, that column is not included in the active spectral certificate, but the active dimension or active rank must be reported explicitly.

The standardized regression is written as
\begin{equation}
\bar{Y}=\bar{\Phi}\bar{K}+E.
\end{equation}
All theoretical results act on these active standardized coordinates. If the full unscreened dictionary contains zero-variance directions, then the minimum eigenvalue of the full Gram matrix is naturally zero. We do not silently hide this degeneracy with regularization; instead, we require reporting active dimensions, active ranks, and spectral quantities.

\section{Data-Quality Diagnostic Certificates}
\label{sec:certificates}

The certificates are defined in the order in which data enter the EDMDc regression. State-space certificates are indirect diagnostics. Lifted- and regression-space certificates are closer to the final estimation problem. Regression isotropy $\creg$ is the directly provable certificate.

\subsection{State-space certificates}

\paragraph{Directional coverage.}
Let
\begin{equation}
d_k=\frac{x_{k+1}-x_k}{\|x_{k+1}-x_k\|_2+\epsilon}
\end{equation}
be a unit displacement direction in the whitened state space. For a set of angular resolutions $\Delta=\{\delta_\ell\}$, let $M_\ell$ be the maximum greedy count of a $\delta_\ell$-separated direction set, and let $M_\ell^\star$ be the reference target count. Define
\begin{equation}
C_{\mathrm{dir}}
=
\min_\ell
\frac{M_\ell}{M_\ell^\star}.
\end{equation}
A larger value indicates more dispersed trajectory directions. This certificate alone does not imply a regression error bound.

\paragraph{Frostman-type non-clustering.}
Let $\tilde{x}_i$ denote a whitened state sample and let $s$ be a reference dimension. For a set of scales $\mathcal R$, define
\begin{equation}
\rho(r)=
\max_i
\frac{
\frac{1}{N}\sum_{j=1}^{N}
\mathbf{1}\{\|\tilde{x}_j-\tilde{x}_i\|_2\le r\}
}{
r^s+\epsilon
}.
\end{equation}
The non-clustering certificate is
\begin{equation}
C_{\mathrm{fr}}
=
\min_{r\in\mathcal R}
\frac{\rho_{\max}}{\rho(r)}.
\end{equation}
This certificate diagnoses local repeated sampling and multiscale clustering. It is an empirical diagnostic and is not equivalent to a strict Frostman measure condition.

\paragraph{State isotropy.}
For the active standardized matrix $\bar X$ of the state sample matrix $X$, define
\begin{equation}
C_{\mathrm{state}}
=
\lambda_{\min}
\left(
\frac{1}{N}\bar{X}^\top\bar{X}
\right).
\end{equation}
In the experiments, \texttt{state\_iso} is this quantity normalized by a fixed display threshold.

\subsection{Lifted- and regression-space certificates}

\paragraph{Lifted-feature isotropy.}
For the active standardized matrix $\bar\Psi$ of the lifted-feature matrix $\Psi=[\psi(x_1)^\top;\ldots;\psi(x_N)^\top]$, define
\begin{equation}
C_{\mathrm{lift}}
=
\lambda_{\min}
\left(
\frac{1}{N}\bar{\Psi}^\top\bar{\Psi}
\right).
\end{equation}
This certificate detects whether the dictionary produces nearly constant columns, zero columns, or highly correlated directions under the sampled distribution.

\paragraph{Regression isotropy.}
For the active standardized matrix $\bar\Phi$ of the concatenated regression matrix $\Phi=[\Psi,U]$, define
\begin{equation}
C_{\mathrm{reg}}(\D_N)
=
\lambda_{\min}
\left(
\frac{1}{N}\bar{\Phi}^\top\bar{\Phi}
\right).
\end{equation}
This certificate acts directly on the final EDMDc least-squares problem. It is not a simple function of \texttt{state\_iso}, \texttt{lift\_iso}, or input variance, because the concatenated matrix can degenerate due to correlations between lifted-feature columns and input columns.

\paragraph{Bottleneck-style composite certificate.}
To report an overall data-quality profile, we use
\begin{equation}
C_{\mathrm{GPE}}
=
\min\{
C_{\mathrm{dir}},
C_{\mathrm{fr}},
C_{\mathrm{rad}},
C_{\mathrm{state}}/\tau_{\mathrm{state}},
C_{\mathrm{lift}}/\tau_{\mathrm{lift}},
C_{\mathrm{reg}}/\tau_{\mathrm{reg}}
\}.
\end{equation}
Here $\tau_{\mathrm{state}}=0.20$, $\tau_{\mathrm{lift}}=0.05$, and $\tau_{\mathrm{reg}}=0.05$ are experimental display thresholds, not theoretical constants. The composite certificate is used to localize bottlenecks and is not a performance guarantee.

\section{Theoretical Guarantees and Boundaries}
\label{sec:theory}

\subsection{Finite-sample regression stability}

Consider the standardized EDMDc regression under a fixed design matrix. For a general finite dictionary, we do not assume that the true system is exactly represented by the dictionary. Define the empirical projection target
\begin{equation}
K_N^\star
=
\arg\min_K
\|\bar{Y}-\bar{\Phi}K\|_F^2,
\end{equation}
and write
\begin{equation}
\bar{Y}=\bar{\Phi}K_N^\star+R_N+E,
\qquad
\bar{\Phi}^\top R_N=0.
\end{equation}
If $\creg(\mathcal D_N)>0$, then
\begin{equation}
\sigma_{\min}(\bar{\Phi})
=
\sqrt{N\creg(\mathcal D_N)}.
\end{equation}
If the noise rows satisfy a conditional Gaussian model $e_k\sim\mathcal{N}(0,\Sigma_e)$ and are mutually independent, then the Fisher-type information matrix is
\begin{equation}
\mathcal I_N
=
\Sigma_e^{-1}\otimes(\bar{\Phi}^{\top}\bar{\Phi}),
\end{equation}
and satisfies
\begin{equation}
\lambda_{\min}(\mathcal I_N)
\ge
\frac{N\creg(\mathcal D_N)}
{\lambda_{\max}(\Sigma_e)}.
\end{equation}
The ridge estimate
\begin{equation}
\widehat K_\lambda
=
(\bar{\Phi}^{\top}\bar{\Phi}+\lambda I)^{-1}\bar{\Phi}^{\top}\bar{Y}
\end{equation}
satisfies
\begin{equation}
\|\widehat K_\lambda-K_N^\star\|_F
\le
\frac{\|\bar{\Phi}^{\top}E\|_F+\lambda\|K_N^\star\|_F}
{N\creg(\mathcal D_N)+\lambda}.
\end{equation}
If the noise rows are conditionally independent, zero-mean, and have covariance bounded by $\sigma^2 I_q$, then the least-squares estimate satisfies
\begin{equation}
\mathbb{E}
\left[
\|\widehat K_0-K_N^\star\|_F^2
\mid
\bar{\Phi}
\right]
\le
\frac{\sigma^2qp}
{N\creg(\mathcal D_N)}.
\end{equation}
This result controls only one-step EDMDc regression stability in active standardized coordinates. It does not control finite-dictionary approximation error, Koopman invariant-subspace error, long-horizon rollout error propagation, or closed-loop control performance.

\subsection{From sampling distributions to empirical regression spectra}

Let $\xi\in\R^p$ be an active regressor sampled from a distribution $\rho$. The theoretical standardization uses population means and scales:
\begin{equation}
\bar{\xi}_\rho
=
D_\rho^{-1}(\xi-\mu_\rho),
\qquad
\mu_\rho=\mathbb{E}_\rho[\xi].
\end{equation}
Define
\begin{equation}
G_\rho
=
\mathbb{E}_\rho[\bar{\xi}_\rho\bar{\xi}_\rho^\top].
\end{equation}
If $\lambda_{\min}(G_\rho)\ge \mu>0$ and independent samples satisfy
\begin{equation}
\left\|
\frac{1}{N}\sum_{k=1}^{N}\bar{\xi}_{\rho,k}\bar{\xi}_{\rho,k}^{\top}
-G_\rho
\right\|_2
\le
\frac{\mu}{2},
\end{equation}
then
\begin{equation}
C_{\mathrm{reg}}^\rho(\mathcal D_N)\ge\frac{\mu}{2}.
\end{equation}
If, additionally, $\|\bar{\xi}_{\rho,k}\|_2\le R$, a matrix Bernstein-type inequality gives that
\begin{equation}
N\gtrsim
\frac{R^2}{\mu}
\log\frac{p}{\delta}
\end{equation}
is sufficient to control the above deviation event with high probability. The sample-standardized \texttt{regression\_iso} used in the experiments is a finite-sample computable substitute for this theoretical object.

\subsection{Identifiability bridge and Schur-complement condition}

Suppose the state distribution has a density lower bound on a compact set $\Omega$, and the active dictionary $\psi_A$ is continuous and linearly independent modulo constants in $L_2(\Omega)$. Then the population Gram matrix of the centered lifted features satisfies
\begin{equation}
A_\rho\succeq \alpha I,
\qquad
\alpha>0.
\end{equation}
If the residual of the input relative to the linear projection onto the lifted features satisfies
\begin{equation}
S_{U|\Psi,\rho}
=
\mathbb{E}
[
(u-\Pi_{\psi_A}u)(u-\Pi_{\psi_A}u)^\top
]
\succeq \beta I,
\end{equation}
then the population Gram matrix of the full regressor satisfies
\begin{equation}
\lambda_{\min}(G_\rho)
\ge
\frac{\min\{\alpha,\beta\}}{(1+\gamma_\rho)^2}.
\end{equation}
This statement shows that a positive regression spectral lower bound requires three classes of conditions to work together: the state distribution must cover the target region, the dictionary must be nondegenerate on that region, and the input must have conditional excitation relative to the lifted features.

\subsection{Counterexample for layer separation}

Consider the two-dimensional sample set
\begin{equation}
\mathcal X=\{(1,0),(-1,0),(0,1),(0,-1)\}.
\end{equation}
Its state covariance is $\Sigma_x=\frac12 I_2$, so the state geometry is favorable. If the quadratic dictionary
\begin{equation}
\psi(x)=[x_1,x_2,x_1^2,x_1x_2,x_2^2]^\top
\end{equation}
is used, then $x_1x_2=0$ for every sample. The column corresponding to $x_1x_2$ in the lifted matrix is identically zero, and both the lifted Gram matrix and the full regression Gram matrix are singular. Thus, state geometry cannot in general imply lifted or regression geometry. This counterexample explains why \texttt{state\_iso}, \texttt{lift\_iso}, and \texttt{regression\_iso} must all be reported.

\subsection{Submodular D-optimal boundary}

Given a finite candidate library $\mathcal V$ of trajectory segments, suppose each candidate segment $i$ contributes a regression block matrix $\Phi_i$. Define the pure regression D-optimal objective
\begin{equation}
F(S)
=
\log\det
\left(
\epsilon I+\sum_{i\in S}\Phi_i^\top\Phi_i
\right).
\end{equation}
This function is monotone and submodular in the set $S$. Therefore, under a fixed candidate library and a fixed budget, the greedy algorithm that selects the largest marginal gain at each step has the classical $(1-1/e)$ approximation guarantee. The full IGPE-DOPT mixed objective also contains state-coverage, directional-novelty, and clustering-penalty terms, and therefore does not inherit this global guarantee.

\section{Certificate-Driven Data Acquisition: IGPE-DOPT}
\label{sec:igpe}

IGPE-DOPT is a certificate-driven sampling example. At each round, the method simulates candidate trajectory segments from a finite library of constant-input segments, appends each candidate segment to the current dataset, computes state, lifted, and regression certificates, and chooses the highest-scoring segment. The algorithm is a myopic set-building procedure, not a globally optimal sampler.

The default scoring terms match the code implementation:
\begin{table}[h]
\centering
\caption{Default IGPE-DOPT scoring terms.}
\label{tab:igpe-weights}
\begin{tabular}{llcl}
\toprule
Term & Code key & Weight & Meaning\\
\midrule
State minimum spectrum & \texttt{state\_min} & 0.75 & State coverage\\
Lifted log-det & \texttt{lift\_logdet} & 0.80 & Lifted D-optimality\\
Lifted effective rank & \texttt{lift\_eff\_rank} & 0.35 & Lifted spectral dispersion\\
Regression minimum spectrum & \texttt{reg\_min} & 0.35 & Regression E-optimality\\
Directional novelty & \texttt{novelty} & 0.30 & New motion directions\\
Input variance & \texttt{u\_var} & 0.10 & Input diversity\\
Clustering penalty & \texttt{cluster} & 0.25 & Non-clustering\\
\bottomrule
\end{tabular}
\end{table}

The comparison methods are split into two categories. \textbf{External/literature-motivated baselines} are methods that represent established or natural alternatives to certificate-driven sampling: \texttt{A-PE} represents adaptive PE-style input probing; \texttt{OID} uses deterministic multisine/chirp probing as an optimal-input-design surrogate; \texttt{GPE-STATE} targets state-space geometric coverage; \texttt{RANDOM} and \texttt{SOBOL} represent undirected random and low-discrepancy sampling; \texttt{STATE-KCENTER} is a state-space space-filling $k$-center coverage baseline whose design intent is to maximize state coverage rather than the regression spectrum; \texttt{LIFT-DOPT} isolates lifted-space log-determinant design; \texttt{REG-DOPT} isolates regression-space log-determinant design; and \texttt{REG-EOPT} is a regression-space E-optimal design baseline whose design intent is to directly improve the minimum regression spectrum.

\textbf{Internal ablation variants} are used only to diagnose the contribution of components inside the proposed scoring rule. These include \texttt{IGPE-NO-DOPT}, \texttt{IGPE-NO-DIR}, \texttt{IGPE-NO-CLUSTER}, and the weight variants \texttt{IGPE-WHALF}, \texttt{IGPE-UNIFORM}, \texttt{IGPE-REG-HEAVY}, and \texttt{IGPE-CLUSTER-HEAVY}. They are not presented as external published baselines; they are controlled variants of the proposed sampler. IGPE-DOPT combines the external diagnostic layers through the mixed certificate score above.

\section{Empirical Validation and Diagnostic Analysis}
\label{sec:experiments}

\subsection{Experimental setup and quality checks}

The experimental systems are Duffing, Van der Pol, and Lorenz. The main budget ablation uses the \texttt{major-budget-ablation} configuration: methods are \texttt{RANDOM}, \texttt{SOBOL}, \texttt{STATE-KCENTER}, \texttt{REG-EOPT}, and \texttt{IGPE-DOPT}; seeds are 0 through 9; budgets are $B\in\{8,12,20,40,80\}$; the segment length is \texttt{l\_seg=12}; the time step is \texttt{dt=0.01}; the prediction horizon is 200 steps; and the control horizon is 100 steps. This configuration produces 750 cases in total.

The expanded baseline comparison uses the separate \texttt{major-revision} configuration: methods are \texttt{RANDOM}, \texttt{SOBOL}, \texttt{STATE-KCENTER}, \texttt{LIFT-DOPT}, \texttt{REG-DOPT}, \texttt{REG-EOPT}, \texttt{A-PE}, \texttt{OID}, \texttt{GPE-STATE}, and \texttt{IGPE-DOPT}; seeds are 0 through 19; the budget is fixed at $B=40$; and all other simulation settings are the same. This comparison is used to assess how certificate diagnostics behave against input-rich, state-covering, lifted-space, and regression-space acquisition baselines.

The quality-check file \url{results/geometry_major_revision_budget/tables/table5_quality_checks.csv} shows that all certificates, spectral quantities, one-step errors, prediction metrics, and control metrics have 750 rows; nonfinite and NaN counts are all zero; and prediction failure and control failure counts are also zero.

\subsection{Scientific question 1: Do certificate layers separate?}

Table~\ref{tab:cert-hierarchy} reports representative compressed results for the three layers: state, lifted, and regression certificates, together with active rank and active dimension. The quantities \texttt{state\_iso}, \texttt{lift\_iso}, and \texttt{regression\_iso} use different display-threshold normalizations and should not be compared as absolute magnitudes across certificate types.

\begin{table}[p]
\centering
\caption{Certificate-layer separation, reported as mean [95\% CI]. The complete table is available at \texttt{results/geometry\_major\_revision\_budget/tables/table6\_v4\_certificate\_hierarchy.md}.}
\label{tab:cert-hierarchy}
\scriptsize
\resizebox{\textwidth}{!}{
\begin{tabular}{llllllll}
\toprule
system & method & n & state\_iso & lift\_iso & regression\_iso & active\_rank & active\_dim\\
\midrule
duffing & RANDOM & 50 & 2.96 [2.55, 3.35] & 0.0829 [0.0373, 0.137] & 0.0815 [0.0363, 0.135] & 10 [10, 10] & 10 [10, 10]\\
duffing & SOBOL & 50 & 3.17 [2.79, 3.5] & 0.0172 [0.00428, 0.0367] & 0.0171 [0.00422, 0.0365] & 9.98 [9.94, 10] & 10 [10, 10]\\
duffing & STATE-KCENTER & 50 & 4.18 [3.9, 4.41] & 0.0154 [0.00539, 0.0277] & 0.0126 [0.00412, 0.0235] & 10 [10, 10] & 10 [10, 10]\\
duffing & REG-EOPT & 50 & 3.27 [2.89, 3.62] & 0.132 [0.0624, 0.21] & 0.131 [0.0613, 0.209] & 10 [10, 10] & 10 [10, 10]\\
duffing & IGPE-DOPT & 50 & 4.16 [3.87, 4.43] & 0.073 [0.024, 0.135] & 0.0711 [0.0235, 0.132] & 10 [10, 10] & 10 [10, 10]\\
lorenz & RANDOM & 50 & 0.549 [0.521, 0.582] & 0.00572 [0.00487, 0.00679] & 0.00546 [0.00456, 0.00657] & 10 [10, 10] & 10 [10, 10]\\
lorenz & SOBOL & 50 & 0.569 [0.545, 0.597] & 0.00584 [0.00489, 0.00682] & 0.00561 [0.00466, 0.00659] & 10 [10, 10] & 10 [10, 10]\\
lorenz & STATE-KCENTER & 50 & 0.598 [0.579, 0.616] & 0.0114 [0.00938, 0.0136] & 0.0108 [0.00864, 0.0129] & 10 [10, 10] & 10 [10, 10]\\
lorenz & REG-EOPT & 50 & 0.611 [0.58, 0.646] & 0.00551 [0.005, 0.00608] & 0.00543 [0.0049, 0.00599] & 10 [10, 10] & 10 [10, 10]\\
lorenz & IGPE-DOPT & 50 & 0.616 [0.573, 0.652] & 0.0086 [0.00671, 0.0106] & 0.00845 [0.00655, 0.0104] & 10 [10, 10] & 10 [10, 10]\\
vdp & RANDOM & 50 & 2.93 [2.47, 3.37] & 0.0635 [0.0339, 0.0977] & 0.0628 [0.0334, 0.0967] & 10 [10, 10] & 10 [10, 10]\\
vdp & SOBOL & 50 & 2.79 [2.3, 3.29] & 0.0722 [0.042, 0.103] & 0.0715 [0.0416, 0.103] & 10 [10, 10] & 10 [10, 10]\\
vdp & STATE-KCENTER & 50 & 3.26 [2.81, 3.67] & 0.0695 [0.0424, 0.0997] & 0.0685 [0.0415, 0.0985] & 9.96 [9.9, 10] & 10 [10, 10]\\
vdp & REG-EOPT & 50 & 2.96 [2.51, 3.41] & 0.0881 [0.0441, 0.13] & 0.0869 [0.043, 0.129] & 9.8 [9.64, 9.92] & 10 [10, 10]\\
vdp & IGPE-DOPT & 50 & 3.97 [3.66, 4.24] & 0.128 [0.0564, 0.206] & 0.125 [0.0546, 0.204] & 10 [10, 10] & 10 [10, 10]\\
\bottomrule
\end{tabular}}
\end{table}

For Duffing, \texttt{STATE-KCENTER} has high \texttt{state\_iso}, but its \texttt{regression\_iso} is lower than those of \texttt{REG-EOPT} and IGPE-DOPT. For Van der Pol, IGPE-DOPT improves both state and regression certificates. For Lorenz, \texttt{STATE-KCENTER} is more favorable for lifted/regression certificates. This shows that the layer relationship is not a monotone chain, but a data-quality profile jointly determined by the system, dictionary, and sampling strategy.

\begin{figure}[p]
\centering
\includegraphics[width=\textwidth]{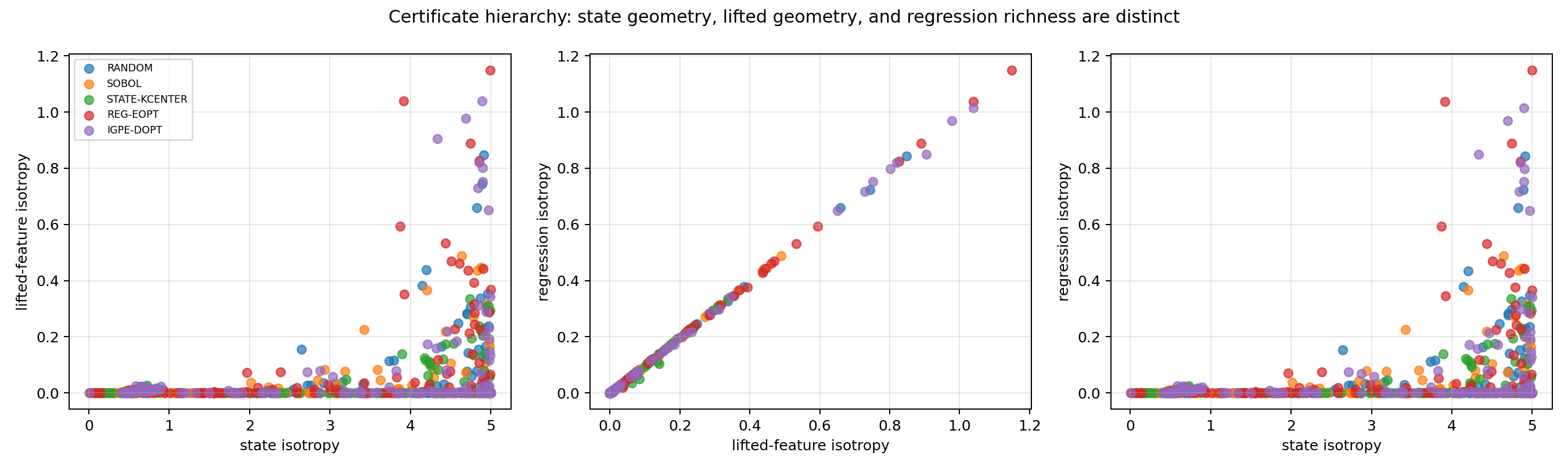}
\caption{Certificate hierarchy and layer separation.}
\label{fig:hierarchy}
\end{figure}

\subsection{Scientific question 2: Does the regression certificate correspond to theoretical quantities?}

Table~\ref{tab:regression-diagnostics} reports the regression certificate, condition number, one-step lift/state RMSE, and active rank. The table does not treat the definitional identity between $\creg$ and $\sigma_{\min}(\bar\Phi)$ as independent empirical evidence; that sanity check is placed in Figure~\ref{fig:a1}.

\begin{table}[p]
\centering
\caption{Expanded baseline comparison at budget $B=40$, reported as mean [95\% CI] over 20 seeds. The underlying summary is available at \texttt{results/geometry\_major\_revision\_baselines/tables/table1\_summary.csv}; the regression-diagnostics table for the design baselines is available at \texttt{results/geometry\_major\_revision\_baselines/tables/table7\_v4\_external\_baselines.md}.}
\label{tab:regression-diagnostics}
\scriptsize
\resizebox{\textwidth}{!}{
\begin{tabular}{llllllll}
\toprule
system & method & n & std\_GPE & regression\_iso & one\_step\_lift & open\_loop & tracking\\
\midrule
duffing & RANDOM & 20 & 0.0565 [0.0299, 0.0888] & 0.0568 [0.0299, 0.0891] & 0.021 [0.0177, 0.0245] & 0.000504 [0.000203, 0.000883] & 0.181 [0.167, 0.193]\\
duffing & SOBOL & 20 & 0.0419 [0.0165, 0.0724] & 0.0419 [0.0165, 0.0724] & 0.0345 [0.0263, 0.0434] & 0.0011 [0.000157, 0.0025] & 0.181 [0.167, 0.193]\\
duffing & STATE-KCENTER & 20 & 0.0066 [0.0034, 0.0111] & 0.0066 [0.0034, 0.0111] & 0.0142 [0.0131, 0.0154] & 0.001 [0.000509, 0.0016] & 0.181 [0.167, 0.193]\\
duffing & LIFT-DOPT & 20 & 0.052 [0.0213, 0.0876] & 0.0628 [0.0234, 0.111] & 0.0537 [0.0386, 0.0692] & 0.0025 [0.000866, 0.0044] & 0.181 [0.167, 0.193]\\
duffing & REG-DOPT & 20 & 0.0775 [0.0422, 0.115] & 0.0847 [0.0455, 0.125] & 0.0453 [0.035, 0.0568] & 0.000504 [0.000238, 0.000813] & 0.181 [0.167, 0.193]\\
duffing & REG-EOPT & 20 & 0.112 [0.0655, 0.156] & 0.127 [0.0696, 0.183] & 0.0242 [0.0177, 0.032] & 0.000965 [0.000321, 0.002] & 0.181 [0.167, 0.193]\\
duffing & A-PE & 20 & 0.026 [0.0101, 0.0489] & 0.0301 [0.0101, 0.0605] & 0.021 [0.0177, 0.0252] & 0.001 [0.000118, 0.0022] & 0.181 [0.167, 0.193]\\
duffing & OID & 20 & 0.0439 [0.027, 0.0608] & 0.0439 [0.027, 0.0608] & 0.0035 [0.0027, 0.0045] & 0.000329 [0.00021, 0.000474] & 0.181 [0.167, 0.193]\\
duffing & GPE-STATE & 20 & 0.0934 [0.0527, 0.142] & 0.112 [0.0596, 0.177] & 0.013 [0.0113, 0.0146] & 0.000459 [0.000159, 0.000856] & 0.181 [0.167, 0.193]\\
duffing & IGPE-DOPT & 20 & 0.0757 [0.0354, 0.121] & 0.0819 [0.0384, 0.129] & 0.053 [0.0395, 0.0686] & 0.309 [0.000752, 0.923] & 0.181 [0.167, 0.194]\\
lorenz & RANDOM & 20 & 0.0053 [0.0045, 0.0061] & 0.0053 [0.0045, 0.0061] & 0.0176 [0.0169, 0.0183] & 12.2 [9.8, 14.3] & 2.54 [0.642, 4.71]\\
lorenz & SOBOL & 20 & 0.0067 [0.0054, 0.0081] & 0.0067 [0.0054, 0.0081] & 0.0173 [0.0168, 0.0179] & 9.3 [7, 11.6] & 3.58 [1.42, 6.06]\\
lorenz & STATE-KCENTER & 20 & 0.0205 [0.0197, 0.0215] & 0.0205 [0.0197, 0.0215] & 0.0368 [0.0358, 0.038] & 7.52 [6.2, 9.23] & 8.5 [5.07, 11.8]\\
lorenz & LIFT-DOPT & 20 & 0.0171 [0.0151, 0.0189] & 0.0171 [0.0151, 0.0189] & 0.0299 [0.0275, 0.0322] & 6.29 [4.85, 7.93] & 2.81 [0.818, 5.54]\\
lorenz & REG-DOPT & 20 & 0.0053 [0.0047, 0.006] & 0.0053 [0.0047, 0.006] & 0.018 [0.0176, 0.0184] & 10.5 [7.98, 13] & 4.77 [2.5, 7.37]\\
lorenz & REG-EOPT & 20 & 0.0069 [0.0059, 0.0079] & 0.0069 [0.0059, 0.0079] & 0.0191 [0.0182, 0.02] & 11.9 [9.54, 14.1] & 3.18 [1.21, 5.47]\\
lorenz & A-PE & 20 & 0.0064 [0.0055, 0.0075] & 0.0064 [0.0055, 0.0075] & 0.0181 [0.0174, 0.0189] & 10.4 [8.12, 12.6] & 3.67 [1.57, 5.99]\\
lorenz & OID & 20 & 0.0044 [0.0041, 0.0049] & 0.0044 [0.0041, 0.0049] & 0.0169 [0.0162, 0.0176] & 12.5 [10.4, 14.5] & 1.45 [0.184, 3.3]\\
lorenz & GPE-STATE & 20 & 0.0068 [0.0059, 0.0079] & 0.0068 [0.0059, 0.0079] & 0.0184 [0.0177, 0.0192] & 12.1 [9.64, 14.2] & 9.15 [6.23, 12]\\
lorenz & IGPE-DOPT & 20 & 0.0102 [0.0085, 0.0119] & 0.0102 [0.0085, 0.0119] & 0.0226 [0.0209, 0.0243] & 8.69 [6.81, 10.8] & 0.787 [0.113, 2.13]\\
vdp & RANDOM & 20 & 0.0497 [0.0374, 0.0623] & 0.0497 [0.0374, 0.0623] & 0.0123 [0.0103, 0.0145] & 0.0529 [0.0086, 0.138] & 0.181 [0.167, 0.193]\\
vdp & SOBOL & 20 & 0.0711 [0.0581, 0.0853] & 0.0711 [0.0581, 0.0853] & 0.0139 [0.0131, 0.0148] & 0.0181 [0.0098, 0.0291] & 0.181 [0.167, 0.193]\\
vdp & STATE-KCENTER & 20 & 0.0804 [0.0599, 0.101] & 0.0804 [0.0599, 0.101] & 0.0096 [0.0087, 0.0104] & 0.0337 [0.0178, 0.0571] & 0.18 [0.166, 0.193]\\
vdp & LIFT-DOPT & 20 & 0.077 [0.0447, 0.114] & 0.0791 [0.0455, 0.116] & 0.0445 [0.0286, 0.062] & 1.93 [0.0993, 4.77] & 0.179 [0.165, 0.192]\\
vdp & REG-DOPT & 20 & 0.0697 [0.0381, 0.106] & 0.0721 [0.0384, 0.113] & 0.0409 [0.0253, 0.0585] & 0.115 [0.0189, 0.288] & 0.181 [0.166, 0.193]\\
vdp & REG-EOPT & 20 & 0.0819 [0.0407, 0.139] & 0.115 [0.0442, 0.214] & 0.0219 [0.0168, 0.0272] & 0.136 [0.0319, 0.324] & 0.18 [0.167, 0.192]\\
vdp & A-PE & 20 & 0.0452 [0.0346, 0.0557] & 0.0452 [0.0346, 0.0557] & 0.0119 [0.0105, 0.0136] & 0.0105 [0.0068, 0.0151] & 0.181 [0.167, 0.193]\\
vdp & OID & 20 & 0.0135 [0.006, 0.0244] & 0.0135 [0.006, 0.0244] & 0.0019 [0.0016, 0.0023] & 0.0355 [0.0198, 0.056] & 0.181 [0.167, 0.193]\\
vdp & GPE-STATE & 20 & 0.0274 [0.012, 0.0442] & 0.0274 [0.012, 0.0442] & 0.0124 [0.0096, 0.0151] & 0.27 [0.0321, 0.612] & 0.18 [0.166, 0.192]\\
vdp & IGPE-DOPT & 20 & 0.108 [0.0633, 0.158] & 0.137 [0.0723, 0.237] & 0.0398 [0.0277, 0.0539] & 4.53 [0.054, 13.4] & 0.187 [0.168, 0.212]\\
\bottomrule
\end{tabular}}
\end{table}

The expanded comparison confirms that the baseline choice changes the diagnostic profile. REG-EOPT and GPE-STATE are favorable on Duffing regression certificates, STATE-KCENTER and LIFT-DOPT are favorable on Lorenz, and IGPE-DOPT has the largest standardized GPE index on Van der Pol. OID gives very small one-step lifted errors on Duffing and Van der Pol, but this does not imply the largest regression certificate. IGPE-DOPT obtains the lowest Lorenz tracking error in this comparison, while it is not best on Duffing or Van der Pol open-loop RMSE. These results are consistent with the theoretical boundary: \texttt{regression\_iso} is the spectral bottleneck for one-step regression stability, but task metrics also depend on dictionary bias, operating region, learned operator spectrum, and closed-loop controllability.

\begin{figure}[p]
\centering
\includegraphics[width=\textwidth]{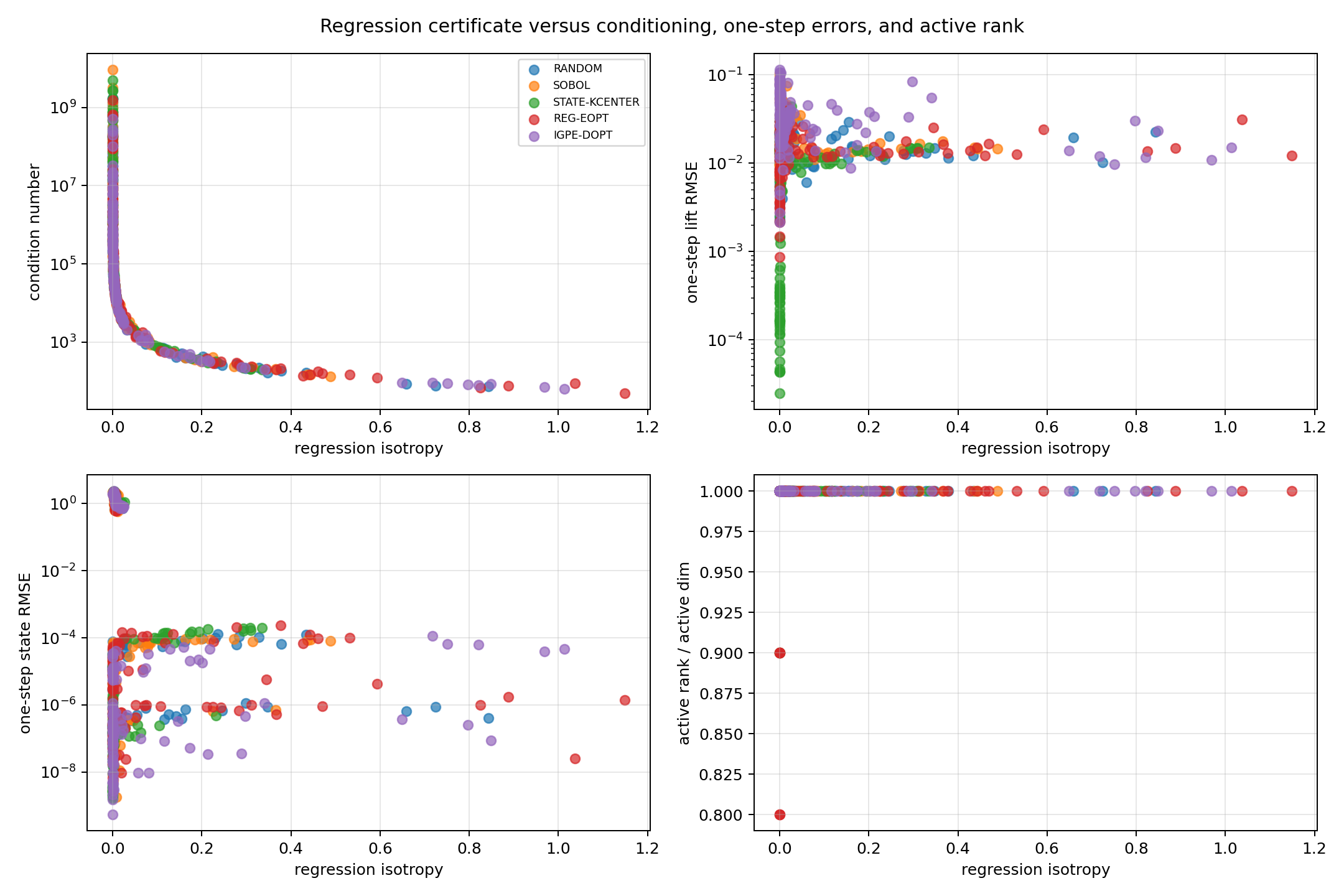}
\caption{Regression certificate and quantities associated with the theory.}
\label{fig:theory-validation}
\end{figure}

\subsection{Scientific question 3: How do budget and sampling distribution affect spectral bottlenecks?}

\begin{figure}[p]
\centering
\includegraphics[width=\textwidth]{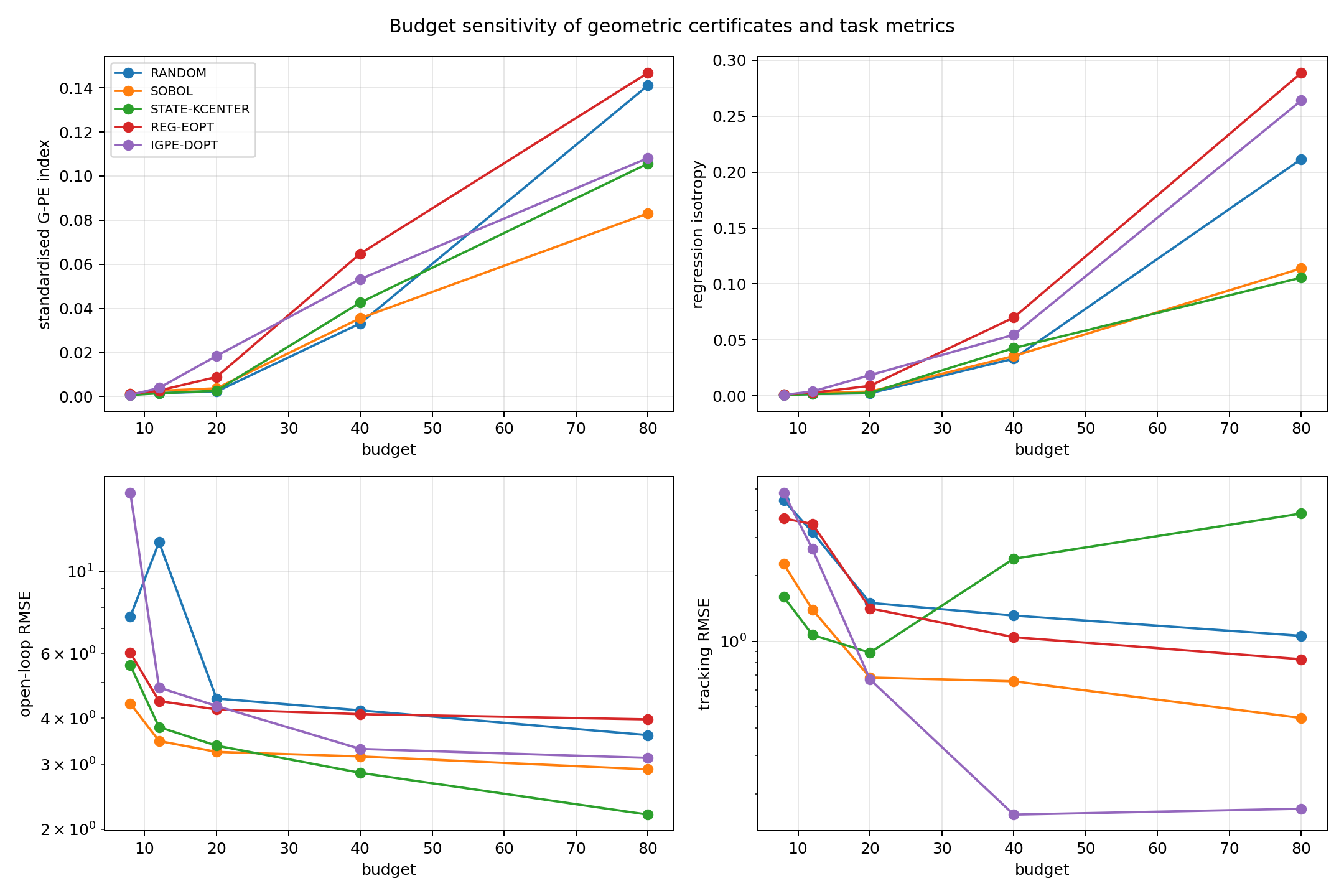}
\caption{Budget sensitivity.}
\label{fig:budget}
\end{figure}

The budget ablation covers $B\in\{8,12,20,40,80\}$ and 750 cases. Figure~\ref{fig:budget} shows that increasing the budget generally improves geometric/regression certificates, but different methods follow different paths. \texttt{REG-EOPT} directly targets the minimum regression spectrum. \texttt{STATE-KCENTER} changes the state distribution. IGPE-DOPT attempts to balance state-, lifted-, and regression-space certificates. This is consistent with the population-to-sample theory: the finite-sample spectral bottleneck is affected both by sample size and by the spectral gap of the population Gram matrix induced by the sampling distribution.

\subsection{Scientific question 4: Are long-horizon prediction and closed-loop tasks monotone?}

Table~\ref{tab:downstream} reports open-loop RMSE, tracking RMSE, failure rates, and paired effect sizes relative to IGPE-DOPT. It shows that the regression certificate cannot be interpreted as a monotone guarantee for long-horizon prediction or closed-loop control.

\begin{table}[p]
\centering
\caption{Downstream prediction/control metrics, reported as mean [95\% CI]. The complete table is available at \texttt{results/geometry\_major\_revision\_budget/tables/table8\_v4\_downstream\_tasks.md}.}
\label{tab:downstream}
\scriptsize
\resizebox{\textwidth}{!}{
\begin{tabular}{llllllll}
\toprule
system & method & open\_loop\_rmse & tracking\_rmse & prediction\_failure\_rate & control\_failure\_rate & dz\_open\_loop\_vs\_IGPE & dz\_tracking\_vs\_IGPE\\
\midrule
duffing & RANDOM & 0.00557 [0.00305, 0.00893] & 0.181 [0.176, 0.186] & 0 & 0 & -0.128 & 0.196\\
duffing & SOBOL & 0.00731 [0.00134, 0.0172] & 0.181 [0.176, 0.186] & 0 & 0 & -0.0619 & -0.0791\\
duffing & STATE-KCENTER & 0.00188 [0.000984, 0.00309] & 0.181 [0.176, 0.186] & 0 & 0 & -0.225 & 0.329\\
duffing & REG-EOPT & 0.0037 [0.00146, 0.00644] & 0.181 [0.176, 0.186] & 0 & 0 & -0.169 & 0.183\\
duffing & IGPE-DOPT & 0.0106 [0.00322, 0.0232] & 0.181 [0.176, 0.186] & 0 & 0 & 0 & 0\\
lorenz & RANDOM & 13.9 [11.6, 17.3] & 6.53 [4.87, 8.33] & 0 & 0 & -0.0536 & 0.196\\
lorenz & SOBOL & 9.82 [8.38, 11.3] & 2.91 [1.73, 4.33] & 0 & 0 & -0.186 & -0.221\\
lorenz & STATE-KCENTER & 10.6 [8.38, 13.6] & 5.53 [3.76, 7.38] & 0 & 0 & -0.157 & 0.0812\\
lorenz & REG-EOPT & 12.2 [11.1, 13.3] & 5.87 [3.97, 7.8] & 0 & 0 & -0.108 & 0.121\\
lorenz & IGPE-DOPT & 15.6 [10.5, 25.2] & 4.73 [2.99, 6.66] & 0 & 0 & 0 & 0\\
vdp & RANDOM & 5.26 [0.768, 13.2] & 0.182 [0.176, 0.19] & 0 & 0 & 0.0549 & 0.0898\\
vdp & SOBOL & 0.471 [0.127, 0.939] & 0.18 [0.176, 0.185] & 0 & 0 & -0.208 & -0.114\\
vdp & STATE-KCENTER & 0.0856 [0.0479, 0.138] & 0.181 [0.176, 0.186] & 0 & 0 & -0.236 & 0.113\\
vdp & REG-EOPT & 1.48 [0.432, 2.73] & 0.208 [0.182, 0.254] & 0 & 0 & -0.134 & 0.184\\
vdp & IGPE-DOPT & 3.54 [0.603, 8.53] & 0.181 [0.176, 0.186] & 0 & 0 & 0 & 0\\
\bottomrule
\end{tabular}}
\end{table}

On Duffing and Van der Pol, the open-loop mean of IGPE-DOPT is affected by a small number of large-error seeds. On Lorenz, the tracking RMSE of IGPE-DOPT is lower than those of \texttt{RANDOM}, \texttt{STATE-KCENTER}, and \texttt{REG-EOPT}, but its open-loop RMSE is not optimal. This is an empirical manifestation of the boundary of the theoretical result: one-step regression stability is not the same as long-horizon rollout stability, nor is it the same as optimality on a closed-loop task.

\begin{figure}[p]
\centering
\includegraphics[width=\textwidth]{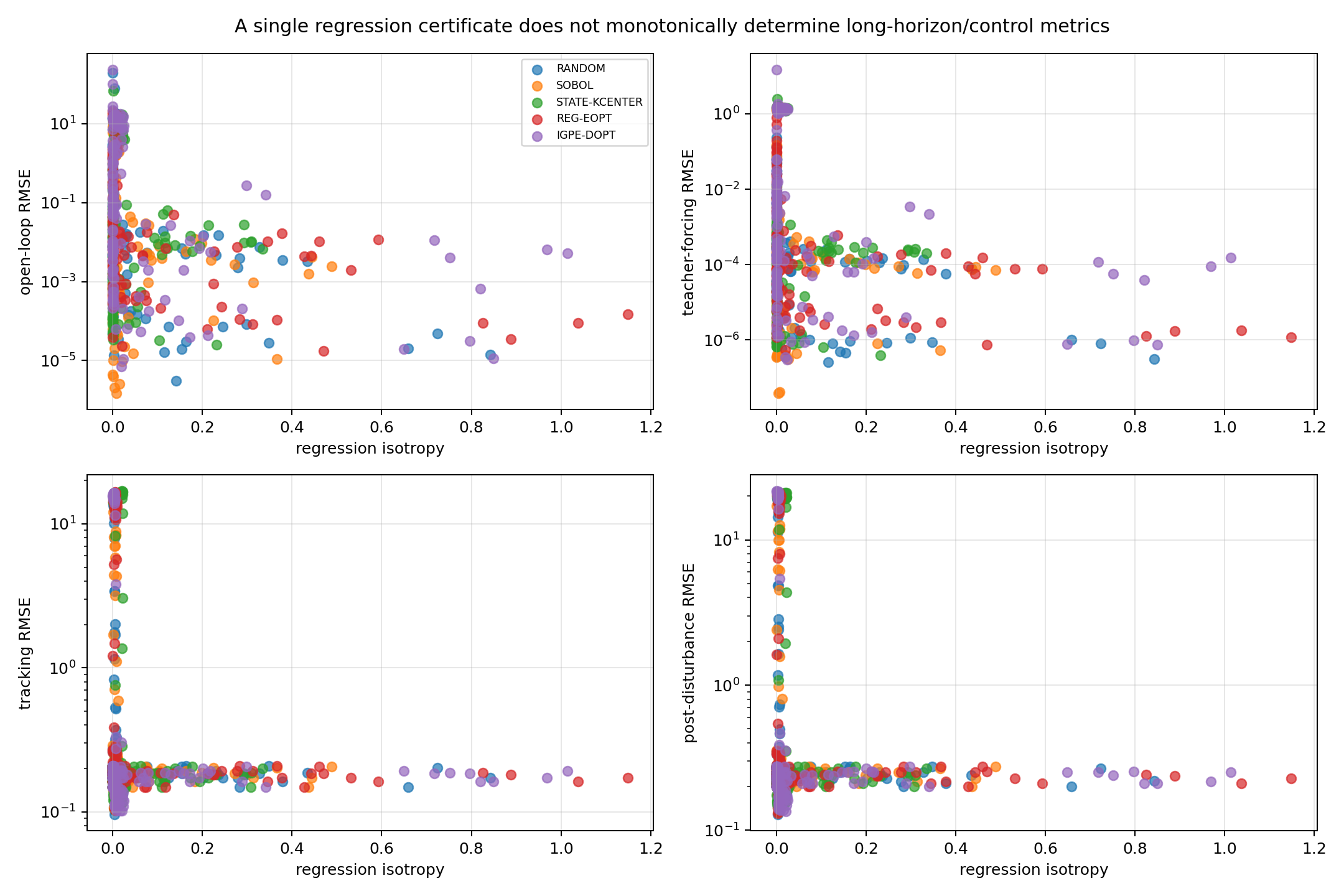}
\caption{Nonmonotonicity in downstream tasks.}
\label{fig:nonmonotonicity}
\end{figure}

\subsection{Scientific question 5: Weight sensitivity}

Table~\ref{tab:weight-sensitivity} comes from the independent \texttt{weight-sensitivity} configuration: three systems, 10 seeds, budget 40, and five IGPE weight variants. The table checks whether the conclusions depend on the default manually tuned weights.

\begin{table}[p]
\centering
\caption{IGPE weight sensitivity, reported as mean values. The complete mean [95\% CI] table is available at \texttt{results/geometry\_weight\_sensitivity/tables/table9\_v4\_weight\_sensitivity.md}.}
\label{tab:weight-sensitivity}
\scriptsize
\setlength{\tabcolsep}{2.5pt}
\begin{tabular}{llrrrr}
\toprule
system & method & \texttt{regression\_iso} & \texttt{one\_step\_lift\_rmse} & \texttt{open\_loop\_rmse} & \texttt{tracking\_rmse}\\
\midrule
duffing & IGPE-DOPT & 0.0698 & 0.0523 & 0.00134 & 0.181\\
duffing & IGPE-WHALF & 0.0812 & 0.0375 & 0.0481 & 0.181\\
duffing & IGPE-UNIFORM & 0.113 & 0.0371 & 0.171 & 0.181\\
duffing & IGPE-REG-HEAVY & 0.149 & 0.0474 & 0.000791 & 0.181\\
duffing & IGPE-CLUSTER-HEAVY & 0.135 & 0.0334 & 0.0302 & 0.181\\
lorenz & IGPE-DOPT & 0.00944 & 0.0216 & 9.72 & 0.119\\
lorenz & IGPE-WHALF & 0.0142 & 0.0254 & 6.83 & 3.03\\
lorenz & IGPE-UNIFORM & 0.00895 & 0.0217 & 11.8 & 2.02\\
lorenz & IGPE-REG-HEAVY & 0.00849 & 0.0211 & 10.1 & 0.132\\
lorenz & IGPE-CLUSTER-HEAVY & 0.0121 & 0.0237 & 11.8 & 4.12\\
vdp & IGPE-DOPT & 0.0844 & 0.0387 & 0.194 & 0.181\\
vdp & IGPE-WHALF & 0.0551 & 0.0330 & 9.10 & 0.179\\
vdp & IGPE-UNIFORM & 0.167 & 0.0302 & 0.465 & 0.182\\
vdp & IGPE-REG-HEAVY & 0.0522 & 0.0404 & 0.116 & 0.181\\
vdp & IGPE-CLUSTER-HEAVY & 0.0665 & 0.0413 & 1.28 & 0.182\\
\bottomrule
\end{tabular}
\end{table}

Changing weights changes the data-quality profile. In Duffing, \texttt{IGPE-REG-HEAVY} increases \texttt{regression\_iso} to 0.149 [0.0772, 0.226]. In Lorenz, \texttt{IGPE-WHALF} increases \texttt{regression\_iso} to 0.0142 [0.0115, 0.0165], but tracking RMSE rises to 3.03 [0.137, 7.3]. This shows that increasing one certificate does not guarantee optimality across all downstream tasks.

\subsection{Computational cost}

No manuscript-scale per-method timing logs were available in the existing result directories, so we ran a lightweight local benchmark on the same script path. The benchmark used Windows 10, Python 3.8.10, an Intel Core Ultra 5 125H CPU with 14 cores and 18 logical processors, and 31.5 GiB RAM. Each entry in Table~\ref{tab:computational-cost} is the mean wall-clock time over three single-system cases (Duffing, Van der Pol, and Lorenz) with seed 0, budget $B=8$, segment length 8, time step 0.01, prediction horizon 40, and control horizon 20. The timings call the internal case runner and exclude file export and figure generation, so they should be interpreted as empirical low-dimensional guidance for relative cost rather than as manuscript-scale runtimes.

\begin{table}[h]
\centering
\caption{Lightweight computational-cost benchmark. Times are seconds per single-system case, reported as mean and range over three systems.}
\label{tab:computational-cost}
\scriptsize
\begin{tabular}{lcc}
\toprule
method & mean wall-clock & range\\
\midrule
RANDOM & 0.151 & 0.145--0.155\\
SOBOL & 0.150 & 0.147--0.156\\
STATE-KCENTER & 0.163 & 0.162--0.164\\
LIFT-DOPT & 0.555 & 0.552--0.561\\
REG-DOPT & 0.561 & 0.558--0.567\\
REG-EOPT & 0.563 & 0.553--0.579\\
A-PE & 0.174 & 0.145--0.206\\
OID & 0.149 & 0.143--0.153\\
GPE-STATE & 0.168 & 0.162--0.172\\
IGPE-DOPT & 0.468 & 0.457--0.476\\
\bottomrule
\end{tabular}
\end{table}

\subsection{Appendix figure for a definitional identity}

\begin{figure}[p]
\centering
\includegraphics[width=0.85\textwidth]{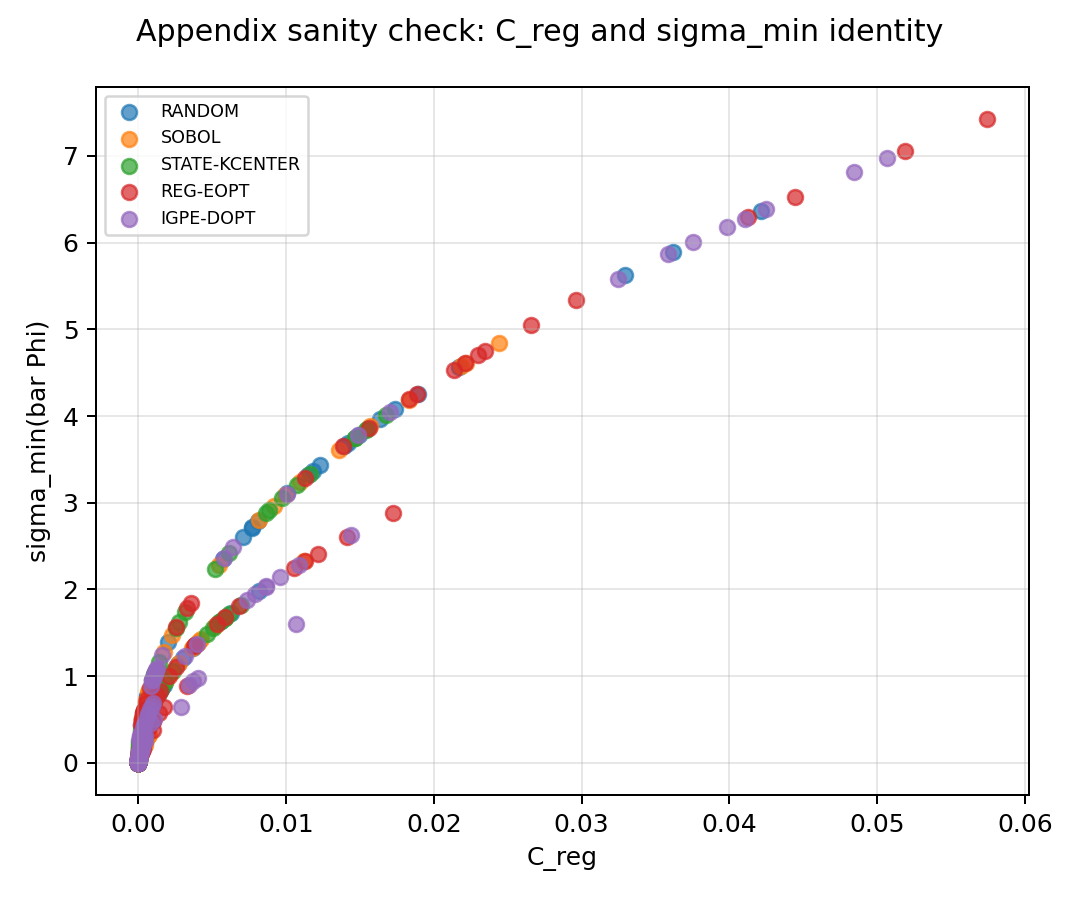}
\caption{Sanity check for the definitional identity between $\creg$ and the minimum singular value.}
\label{fig:a1}
\end{figure}

Figure~\ref{fig:a1} shows the sanity check between $\creg$ and $\sigma_{\min}(\bar\Phi)$. Because they satisfy the definitional identity $\sigma_{\min}(\bar\Phi)=\sqrt{N\creg}$, the figure is used only as an implementation check and not as an independent empirical finding.

\section{Discussion and Limitations}
\label{sec:discussion}

The claim of this paper is diagnostic. Data quality for Koopman/EDMDc should not ask only whether the input is PE; it should ask whether the population Gram matrix and finite-sample Gram matrix of the active standardized regressor have sufficient spectral gaps. State coverage, non-clustering, and directional diversity explain the sampling-distribution layer. Lifted and regression isotropy explain the dictionary layer and the final regression layer.

The paper has four main boundaries. First, we do not prove that directional coverage or Frostman non-clustering implies a lower bound on $\creg$ for general dictionaries and general dynamical systems. Second, IGPE-DOPT is a mixed heuristic sampler; the full objective has no global optimality or approximate-optimality guarantee. Third, the regression theory controls only one-step EDMDc estimation stability, not long-horizon prediction, closed-loop control, or Koopman invariant-subspace recovery. Fourth, all spectral certificates depend on a fixed dictionary and active-column rules; therefore, reports must jointly provide \texttt{active\_dim}, \texttt{active\_rank}, and spectral quantities.

These boundaries are not incidental caveats; they are part of the argument. Table~\ref{tab:downstream} and Figure~\ref{fig:nonmonotonicity} show that downstream task metrics are not monotone in a single regression certificate. Table~\ref{tab:weight-sensitivity} shows that weight changes alter the certificate profile, but do not produce a single task-optimal setting. The role of the certificates is to localize data-quality bottlenecks, not to replace task-level validation.

\paragraph{Practitioner guide.}
In the low-dimensional experiments of this paper, the certificates are most useful as a triage workflow. If active rank is deficient, the first issue is structural: the chosen dictionary, sampling region, or input protocol has produced inactive or redundant regression columns, and adding a small amount of regularization should not be reported as if the design were informative. If active rank is full but \texttt{regression\_iso} or $\creg$ is small relative to other methods under the same system, dictionary, and budget, the bottleneck is the final EDMDc regression spectrum and a regression-space design such as \texttt{REG-EOPT} or the regression term in IGPE-DOPT is the natural repair. If state certificates are poor while lifted/regression spectra are also poor, state-space sampling should be repaired first through coverage-oriented strategies such as \texttt{STATE-KCENTER}, \texttt{GPE-STATE}, or the state terms in IGPE-DOPT. If state coverage is high but lifted or regression spectra remain poor, the likely bottleneck is dictionary degeneracy or lack of conditional input excitation relative to lifted features. In all cases, the thresholds are empirical rules for the present low-dimensional settings: compare certificates only after fixing the system, dictionary, standardization rule, budget, and active-column threshold, and then validate the selected design on the downstream prediction or control task.

\paragraph{High-dimensional scalability.}
The current implementation is intended as a diagnostic prototype, not a high-dimensional optimized solver. Spectral certificates require covariance or Gram-matrix computations whose cost grows with the active lifted dimension. Candidate trajectory scoring multiplies this cost by the number of proposed segments and acquisition rounds. Non-clustering and pairwise coverage diagnostics can become quadratic in the number of samples unless approximate nearest-neighbor, sketching, or minibatch estimators are introduced. Polynomial dictionaries further increase the lifted dimension combinatorially with state dimension and degree, so high-dimensional systems require either structured dictionaries, randomized features, sparse libraries, or incremental spectral updates.

\paragraph{Code and artifact availability.}
The code, scripts, figures, tables, and LaTeX source used for this revision are prepared as supplementary material in the artifact package accompanying the manuscript. No external repository link is cited in the manuscript; if a public archive is released later, this sentence can be replaced with the archived link and identifier.

\paragraph{Broader impact and ethics.}
The proposed certificates are diagnostics for data quality and regression identifiability. They are not safety guarantees. In safety-critical control, a favorable certificate profile must still be combined with constraints, robustness analysis, independent validation, failure-mode testing, and human engineering review before deployment.

\section{Conclusion}
\label{sec:conclusion}

This paper reformulates Koopman/EDMDc data acquisition as a problem of data-quality diagnostics and regression identifiability. Input PE, state coverage, lifted-feature richness, and final EDMDc regression identifiability are different layers and cannot replace one another.

We defined multilayer diagnostic certificates, proved that $\creg$ has direct meaning for one-step EDMDc regression stability, and established a theoretical chain from sampling distributions to empirical regression spectral lower bounds. The state-geometry counterexample and the Schur-complement condition show that the concatenated regression matrix must be inspected directly; reporting only state coverage, lifted-feature rank, or input variance is insufficient.

The experiments show that certificate layers do separate, that budget and weights alter the data-quality profile, and that long-horizon prediction and closed-loop control do not vary monotonically with a single regression certificate. IGPE-DOPT demonstrates the feasibility of certificate-driven sampling, but the contribution of the paper is not a universally optimal sampler. The contribution is a computable, interpretable, ablatable diagnostic framework with explicit theoretical boundaries.

\appendix

\section{Reproduction Commands}

\begin{lstlisting}[language=bash,basicstyle=\ttfamily\small]
python code/test_improved_gpe.py

python code/improved_gpe_experiments.py major-revision \
  --output-dir results/geometry_major_revision_baselines

python code/improved_gpe_experiments.py major-budget-ablation \
  --output-dir results/geometry_major_revision_budget

python code/improved_gpe_experiments.py degree-ablation \
  --output-dir results/geometry_degree_ablation

python code/improved_gpe_experiments.py component-ablation \
  --output-dir results/geometry_component_ablation

python code/improved_gpe_experiments.py weight-sensitivity \
  --output-dir results/geometry_weight_sensitivity
\end{lstlisting}

\section{Code Metrics and Table/Figure Outputs}

\begin{table}[h]
\centering
\caption{Mapping between manuscript objects and code columns.}
\label{tab:code-columns}
\begin{tabular}{ll}
\toprule
Manuscript object & Code column\\
\midrule
$\creg$ & \texttt{regression\_cov\_z\_min}\\
$\sigma_{\min}(\bar{\Phi})$ & \texttt{sigma\_min\_bar\_phi}\\
Regression log-det & \texttt{regression\_logdet}\\
Active rank/dimension & \texttt{active\_rank}, \texttt{active\_dim}\\
One-step regression error & \texttt{one\_step\_lift\_rmse}, \texttt{one\_step\_state\_rmse}\\
Composite bottleneck certificate & \texttt{std\_gpe\_index}\\
\bottomrule
\end{tabular}
\end{table}

\begin{table}[h]
\centering
\caption{Table and figure output files.}
\label{tab:artifact-map}
\scriptsize
\begin{tabular}{p{0.35\textwidth}p{0.58\textwidth}}
\toprule
Table or figure & File\\
\midrule
Table 6 certificate-layer separation & \url{results/geometry_major_revision_budget/tables/table6_v4_certificate_hierarchy.csv}, \url{results/geometry_major_revision_budget/tables/table6_v4_certificate_hierarchy.md}\\
Table 7 expanded baseline comparison & \url{results/geometry_major_revision_baselines/tables/table1_summary.csv}, \url{results/geometry_major_revision_baselines/tables/table7_v4_external_baselines.md}\\
Table 8 downstream prediction/control metrics & \url{results/geometry_major_revision_budget/tables/table8_v4_downstream_tasks.csv}, \url{results/geometry_major_revision_budget/tables/table8_v4_downstream_tasks.md}\\
Table 9 weight sensitivity & \url{results/geometry_weight_sensitivity/tables/table9_v4_weight_sensitivity.csv}, \url{results/geometry_weight_sensitivity/tables/table9_v4_weight_sensitivity.md}\\
Figure 6 budget sensitivity & \url{results/geometry_major_revision_budget/figures/figure6_budget_sensitivity.png}\\
Figure 9 certificate-layer separation & \url{results/geometry_major_revision_budget/figures/figure9_certificate_hierarchy.png}\\
Figure 10 regression-theory validation & \url{results/geometry_major_revision_budget/figures/figure10_regression_theory_validation.png}\\
Figure 11 downstream task nonmonotonicity & \url{results/geometry_major_revision_budget/figures/figure11_task_nonmonotonicity.png}\\
Figure A1 sanity check between $\creg$ and minimum singular value & \url{results/geometry_major_revision_budget/figures/figureA1_creg_sigma_min_sanity.png}\\
\bottomrule
\end{tabular}
\end{table}

\bibliographystyle{plainnat}
\bibliography{references}

\end{document}